\documentclass[11pt]{amsart}
\usepackage[sorted-cites]{amsrefs}
\usepackage{color}
\usepackage{amscd}

\newfont{\fnt}{cmr10 scaled 550}

\newtheorem{theorem}{Theorem}

\newtheorem{lemma}{Lemma}

\theoremstyle{remark}

\font\strange=msbm10

\renewcommand{\epsilon}{\varepsilon}

\renewcommand{\Sigma}{\varSigma}
\newcommand{\R}{{{\mathchoice  {\hbox{$\textstyle{\text{\strange R}}$}}
{\hbox{$\textstyle{\text{\strange R}}$}}
{\hbox{$\scriptstyle  N\kern-0.3em  R$}}  
{\hbox{$\scriptscriptstyle  R\kern-0.2em  R$}}}}}

\newcommand{\Z}{{{\mathchoice  {\hbox{$\textstyle{\text{\strange Z}}$}}
{\hbox{$\textstyle{\text{\strange Z}}$}}
{\hbox{$\scriptstyle  Z\kern-0.3em  Z$}}
{\hbox{$\scriptscriptstyle  Z\kern-0.2em  Z$}}}}}

\newcommand{\N}{{{\mathchoice  {\hbox{$\textstyle{\text{\strange N}}$}}
{\hbox{$\textstyle{\text{\strange N}}$}}
{\hbox{$\scriptstyle  N\kern-0.3em  N$}}
{\hbox{$\scriptscriptstyle  N\kern-0.2em  N$}}}}}

\renewcommand{\phi}{\varphi}
\usepackage{amsmath,amsthm,amscd}

\begin{document}
 \begin{abstract}
In 2005, B\"ottcher and Wenzel raised the conjecture that if $X,Y$ are real square matrices, then $||XY-YX||^2\leq 2||X||^2||Y||^2$, where $||\cdot||$ is the Frobenius norm.  Various proofs of this conjecture were found in
the last few years by several authors. We here give another
proof. This proof is highly conceptual and requires minimal
computation. We also briefly discuss related inequalities,
in particular, the classical Chern-do Camo-Kobayashi inequality.
\end{abstract}

\title{Remarks on the B\"ottcher-Wenzel inequality}

\date{May 31, 2011}


\author[Zhiqin Lu]{Zhiqin Lu}
\address{Department of Mathematics, University of California, Irvine, CA 92697.}
\email[Zhiqin Lu]{zlu@uci.edu}

\thanks{The 
 author is partially supported by  NSF award DMS-0904653.}
 
 \newcommand{\M}{\mathcal M}

\maketitle
\section{The proof}

The B\"ottcher-Wenzel inequality was proved in~\cite{bw} for
$2 \times 2$ matrices and by L\'{a}szl\'{o}~\cite{la} for
$3 \times 3$ matrices. 
Different proofs of the
 full version of the conjecture were obtained 
by Lu~\cite{lu}, Vong and Jin~\cite{vj-1},  B\"ottcher-Wenzel~\cite{bw-2}, and Audenaert~\cite{au}. The complex matrix case was treated in 
B\"ottcher-Wenzel~\cite{bw-2} and Wenzel~\cite{wenzel-2}.  A convenient observation that links the complex case to the real one  can be found  in Cheng-Vong-Wenzel~\cite{cvw}*{page 296}. 
A useful observation was obtained in Audenaert~\cite{aw} for further generalizations. 

In this section, we give a new proof. Let $[X,Y]$ denote the
commutator of $X$ and $Y$ and let $\| \cdot\|$ stand for the Frobenius
norm ($=$ Hilbert-Schmidt norm).

\begin{theorem} \label{thm1}
Let $X,Y$ be real $n\times n$ matrices. Then
\[
||[X,Y]||^2\leq 2||X||^2\cdot||Y||^2.
\]
\end{theorem}

In accordance with \cite{lu}*{Lemma 3}, we made the following definition. 
Let $V={\mathfrak g}{\mathfrak l}(n,\mathbb R)$ and define a linear map $T=T_X$ by
\[
T: V\to V,\quad Y\mapsto [X^T,[X,Y]],
\]
where $X^T$ is the transpose of $X$.

Let
\[
\Lambda=\begin{pmatrix}
s_1\\&\ddots\\&&s_n
\end{pmatrix}
\]
and  assume that
\begin{equation}\label{2}
s_1^2+\cdots+s_n^2=1.
\end{equation}

Put
\[
\Lambda_1=\begin{pmatrix}
0&\Lambda\\
0&0
\end{pmatrix}.
\]
Let  $\tilde T_{\Lambda_1}$ be the restriction of $T_{\Lambda_1}$ on  $V={\mathfrak g}{\mathfrak l}(n,\mathbb R)\oplus{\mathfrak g}{\mathfrak l}(n,\mathbb R)\subset {\mathfrak g}{\mathfrak l}(2n,\mathbb R)$. Let finally
\begin{equation}\label{1}
A=\begin{pmatrix}
C\\
&B\end{pmatrix}.
\end{equation}
Then we have
\begin{equation}\label{3}
[\Lambda_1, A]=\begin{pmatrix} 0 &\Lambda B-C\Lambda\\
0&0
\end{pmatrix},
\end{equation}
and 
\[
\tilde T_{\Lambda_1}(A)= \begin{pmatrix} -(\Lambda B-C\Lambda)\Lambda\\&\Lambda  (\Lambda B-C\Lambda)
\end{pmatrix}.
\]

For the rest of the paper, we make the following generic condition: all $s_i$ are distinct and  nonzero, $s_1^2>s_2^2>\cdots >s_n^2$, and all $s_i^2+s_j^2$ are distinct.

\begin{lemma}\label{lem2}
The eigenvalues of $\tilde T_{\Lambda_1}$ must either be $0$ or of the form $s_i^2+s_j^2$. Let $E_{ij}$ be the matrices whose only nonzero entry $1$ is the $(i,j)$-th entry. Then
\begin{enumerate}
\item The eigenspace of the eigenvalue $2s_i^2$ is spanned by $(B,C)=(E_{ii},-E_{ii})$ for $1\leq i\leq n$;
\item The eigenspace of the eigenvalue $s_i^2+s_j^2$ for $i\neq j$ is spanned by $(B,C)=(E_{ij}, -\frac{s_j}{s_i} E_{ij})$ and $(B,C)=(E_{ji},-\frac{s_j}{s_i} E_{ji})$ for $1\leq i\neq j\leq n$;
\item The eigenspace of the eigenvalue  $0$ is spanned by $(B,C)=(E_{ij},\frac{s_i}{s_j} E_{ij})$ for $1\leq i,j\leq n$.
\end{enumerate}
In particular,
the maximum eigenvalue of $\tilde T_{\Lambda_1}$  is $2s_1^2$, of  multiplicity $1$, and the second largest eigenvalue of $\tilde T_{\Lambda_1}$ is  $s_1^2+s_2^2$,  of multiplicity $2$.
\end{lemma}

{\bf Proof.} Let $A$ in~\eqref{1} be an eigenvector of  the  eigenvalue $\lambda$ of $\tilde T_{\Lambda_1}$. Then we have
\[
-(\Lambda B-C\Lambda)\Lambda=\lambda C,\quad
\Lambda  (\Lambda B-C\Lambda)=\lambda B.
\]
Assuming that $b_{ij}, c_{ij}$ are the entries of $B,C$,  respectively, we have
\begin{equation}\label{6}
-s_is_j b_{ij}=(\lambda-s_j^2) c_{ij},\quad
-s_is_jc_{ij}=(\lambda-s_i^2) b_{ij}
\end{equation}
for $1\leq i,j\leq n$.
From the above equations, we conclude that the eigenvalues of $\tilde T_{\Lambda_1}$ must be the solutions of the equations
\[
s_i^2s_j^2b_{ij}c_{ij}=(\lambda-s_i^2)(\lambda-s_j^2)b_{ij}c_{ij}
\]
and hence be either $0$ or  $s_i^2+s_j^2$. Moreover, for fixed $(i,j)$ and the  fixed eigenvalue $s_i^2+s_j^2$, we have $b_{rs}c_{rs}= 0$  except $(r,s)=(i,j)$ or $(j,i)$. Using this observation, we  find all the eigenvectors of the operator $\tilde T_{\Lambda_1}$.  

\qed

{\bf Proof of the Theroem~\ref{thm1}.} 
Following~\cite{bw}, we work with  the  singular value decomposition. Let $||X||=1$ and 
let
\begin{equation}\label{77}
X=Q_1\Lambda Q_2
\end{equation}
be the singular decomposition of $X$, where $Q_1,Q_2$ are orthogonal matrices and $\Lambda$ is a diagonal matrix. Let
\[
B=Q_2YQ_2^{-1},\quad C=Q_1^{-1}YQ_1.
\]
Then we have
\begin{equation}\label{4}
||[X,Y]||^2=||\Lambda B-C\Lambda||^2.
\end{equation}

For fixed $X$, let $Y$ be a matrix with  unit norm such that $||[X,Y]||$ is maximized. Then we have
\[
T_X(Y)=\lambda' Y
\]
by the  method of Lagrange multipliers. By~\cite{lu}*{Proposition 5} (see also ~\cite{cvw}*{Proposition 2.4}), $[X^T,Y^T]$ is also an eigenvector of $\lambda'$ and it is linearly independent to $Y$. 
Let 
\[
Z=\alpha Y+\beta[X^T,Y^T]
\]
be a linear combination of $Y, [X^T,Y^T]$ such that 
\[
A=\begin{pmatrix} Q_1^{-1}ZQ_1\\&Q_2ZQ_2^{-1}\end{pmatrix}
\]
is orthogonal to the first eigenspace  of $\tilde T_{\Lambda_1}$. Since the space of all such  $A$ is $2$-dimensional,   the  linear combination always exists. By~\eqref{3},~\eqref{4}, we have
\[
||[X,Y]||^2=||[\Lambda_1,A]||^2=\langle A,\tilde T_{\Lambda_1}(A)\rangle.
\]
Using  Lemma~\ref{lem2}, we have
\[
||[X,Y]||^2\leq (s_1^2+s_2^2)||A||^2\leq ||B||^2+||C||^2=2||Y||^2,
\]
and the theorem is proved.

\qed

\section{Additional Remarks}
Here are some remarks on further  generalizations of the B\"ottcher-Wenzel inequality.
We first prove the following result.
\begin{theorem}\label{lu}
Let $X,Y$ be $n\times n$ matrices.
Let $X$ be a diagonal matrix and let $||Y||_{\infty} = \max_{i\neq j} (|y_{ij}|)$, where $(y_{ij})$  are the entries of $Y$. 
Then we have
\begin{equation}\label{lu-2}
||[X,Y]||^2\leq ||X||^2 \cdot (||Y||^2+2||Y||^2_\infty).
\end{equation}
\end{theorem}

{\bf Proof.}
In ~\cite{lu}*{pp 1293, Remark 1}, the theorem was proved  assuming that $Y$ is symmetric.  That is, for any real numbers $\lambda_1,\cdots,\lambda_n$, we have
\[
2\sum_{i<j}(\lambda_i-\lambda_j)^2 y_{ij}^2\leq (\sum_j\lambda_j^2)\cdot (2\sum_{i> j} y_{ij}^2+2\max_{i>j} (y_{ij})^2).
\]
 This implies that, for  strictly upper triangular matrix $Y_1$,
\[
||[X,Y_1]||^2\leq ||X||^2\cdot(||Y_1||^2+||Y_1||^2_\infty).
\]
By using the same argument, the above inequality is also  true for strictly lower triangular matrices. Let $Y=Y_0+Y_1+Y_2$, where $Y_0, Y_1, Y_2$ are the diagonal part, the strictly upper triangular part, and the strictly lower triangular parts of $Y$, respectively.  Then we have
\begin{align*}
&
||[X,Y]||^2=||[X,Y_1]||^2+||[X,Y_2]||^2\\
&\leq ||X||^2\cdot(||Y_1||^2+||Y_1||^2_\infty+||Y_2||^2+||Y_2||^2_\infty)\\&
\leq 
 ||X||^2 \cdot (||Y||^2+2||Y||^2_\infty),
 \end{align*}
 and the theorem is proved.
 
 \qed
 
 Chern-do Carmo-Kobayashi~\cite{chern-d-k} already had the B\"ottcher-Wenzel inequality when one of the matrices is symmetric. Their proof  actually yields  the following result.
\begin{theorem}\label{thm2}
Let $X$ be a symmetric $n\times n$  matrix with $\lambda_1$ being the largest eigenvalue and $\lambda_n$ being the smallest eigenvalue. Let $Y$ be an $n\times n$ matrix. Then
\[
||[X,Y]||^2\leq (\lambda_1-\lambda_n)^2||Y||^2.
\]
\end{theorem}

{\bf Proof.} Without loss of generality, we may assume that $X$ is a diagonal matrix. Thus we have
\[
||[X,Y]||^2=\sum_{i,j}(\lambda_i-\lambda_j)^2y_{ij}^2\leq  (\lambda_1-\lambda_n)^2||Y||^2,
\]
where $(y_{ij})$ are the entries of $Y$.

\qed

Let $X$ be a real $n\times n$ matrix and let $s_1,\cdots, s_n$ be the singular values of $X$. The $(2,(2))$-Ky Fan norm of $X$ is defined as 
\[
||X||_{2,(2)}=\sqrt{s_1^2+s_2^2}.
\]
Obviously, we have $||X||_{2,(2)}\leq ||X||$. From the proof of Theorem~\ref{thm1}, we actually have 
\begin{equation}\label{wenzel}
||[X,Y]||^2\leq 2||X||_{2,(2)}^2||Y||^2,
\end{equation}
which is a generalization of the B\"ottcher-Wenzel inequality. This strengthened inequality  was first proved by
 Wenzel~\cite{wenzel-2}.

Evidently, we have
\[
(\lambda_1-\lambda_n)^2\leq 2\max _{i\neq j} \,(\lambda_i^2+\lambda_j^2).
\]
Therefore, if $X$ is a symmetric matrix, then the Chern-do Carmo-Kobayashi inequality is sharper than the Wenzel inequality~\eqref{wenzel}. 
On the other hand, in a lot of cases inequality~\eqref{lu-2}  is sharper than both the Chern-do Carmo-Kobayashi and the Wenzel's inequalities because the $\infty$ norm is usually much smaller. 
 We wish to obtain a common generalization of the above three inequalities. Such a result would provide a common generalization of the B\"ottcher-Wenzel inequality and the Normal Scalar Curvature inequality proved in~\cite{lu} and~\cite{ge-tang}.

Finally, a generalization of Theorem~\ref{lu} may exist in connection with the $p$ Schatten norms considered in~\cite{cvw,wenzel-2, aw}.\\

{\bf Acknowledgement.} The author deeply thanks Professors A. B\"ottcher and D. Wenzel for their many  useful   comments without which the paper would not be  in its current form.


\begin{bibdiv}
\begin{biblist} 

\bib{au}{article}{
   author={Audenaert, K. M. R.},
   title={Variance bounds, with an application to norm bounds for
   commutators},
   journal={Linear Algebra Appl.},
   volume={432},
   date={2010},
   number={5},
   pages={1126--1143},
   issn={0024-3795},
   review={\MR{2577614 (2011b:15047)}},
   doi={10.1016/j.laa.2009.10.022},
}


\bib{bw}{article}{
   author={B{\"o}ttcher, A.},
   author={Wenzel, D.},
   title={How big can the commutator of two matrices be and how big is it
   typically?},
   journal={Linear Algebra Appl.},
   volume={403},
   date={2005},
   pages={216--228},
   issn={0024-3795},
   review={\MR{2140283 (2006f:15016)}},
   doi={10.1016/j.laa.2005.02.012},
}


\bib{bw-2}{article}{
   author={B{\"o}ttcher, A.},
   author={Wenzel, D.},
   title={The Frobenius norm and the commutator},
   journal={Linear Algebra Appl.},
   volume={429},
   date={2008},
   number={8-9},
   pages={1864--1885},
   issn={0024-3795},
   review={\MR{2446625 (2009j:15089)}},
   doi={10.1016/j.laa.2008.05.020},
}

\bib{cvw}{article}{
   author={Cheng, C-M},
   author={Vong, S-W},
   author={Wenzel, D.},
   title={Commutators with maximal Frobenius norm},
   journal={Linear Algebra Appl.},
   volume={432},
   date={2010},
   number={1},
   pages={292--306},
   issn={0024-3795},
   review={\MR{2566477 (2010m:15031)}},
   doi={10.1016/j.laa.2009.08.008},
}


\bib{chern-d-k}{article}{
   author={Chern, S. S.},
   author={do Carmo, M.},
   author={Kobayashi, S.},
   title={Minimal submanifolds of a sphere with second fundamental form of
   constant length},
   conference={
      title={Functional Analysis and Related Fields (Proc. Conf. for M.
      Stone, Univ. Chicago, Chicago, Ill., 1968)},
   },
   book={
      publisher={Springer},
      place={New York},
   },
   date={1970},
   pages={59--75},
   review={\MR{0273546 (42 \#8424)}},
}


\bib{flc}{article}{
   author={Fong, K-S},
   author={Cheng, C-M},
   author={Lok, I-K},
   title={Another unitarily invariant norm attaining the minimum norm bound
   for commutators},
   journal={Linear Algebra Appl.},
   volume={433},
   date={2010},
   number={11-12},
   pages={1793--1797},
   issn={0024-3795},
   review={\MR{2736098}},
   doi={10.1016/j.laa.2010.06.037},
}

\bib{ge-tang}{article}{
   author={Ge, J.},
   author={Tang, Z.},
   title={A proof of the DDVV conjecture and its equality case},
   journal={Pacific J. Math.},
   volume={237},
   date={2008},
   number={1},
   pages={87--95},
   issn={0030-8730},
   review={\MR{2415209 (2009d:53080)}},
   doi={10.2140/pjm.2008.237.87},
}

\bib{la}{article}{
   author={L{\'a}szl{\'o}, L.},
   title={Proof of B\"ottcher and Wenzel's conjecture on commutator norms
   for 3-by-3 matrices},
   journal={Linear Algebra Appl.},
   volume={422},
   date={2007},
   number={2-3},
   pages={659--663},
   issn={0024-3795},
   review={\MR{2305148 (2008a:15053)}},
   doi={10.1016/j.laa.2006.11.021},
}


\bib{lu}{article}{
author={Lu,Z.},
title={Normal Scalar Curvature Conjecture and its applications},
journal={Journal of Functional Analysis},
volume={261},
date={2011},
pages={1284--1308},
 doi={doi:10.1016/j.jfa.2011.05.002},}


\bib{vj-1}{article}{
   author={Vong, S-W},
   author={Jin, X-Q},
   title={Proof of B\"ottcher and Wenzel's conjecture},
   journal={Oper. Matrices},
   volume={2},
   date={2008},
   number={3},
   pages={435--442},
   issn={1846-3886},
   review={\MR{2440678 (2009f:15033)}},
}

\bib{wenzel-2}{article}{
   author={Wenzel, D.},
   title={Dominating the commutator},
   conference={
      title={Topics in operator theory. Volume 1. Operators, matrices and
      analytic functions},
   },
   book={
      series={Oper. Theory Adv. Appl.},
      volume={202},
      publisher={Birkh\"auser Verlag},
      place={Basel},
   },
   date={2010},
   pages={579--600},
   review={\MR{2723302 (2011j:15034)}},
}
		

\bib{aw}{article}{
   author={Wenzel, D.},
   author={Audenaert, K. M. R.},
   title={Impressions of convexity: an illustration for commutator bounds},
   journal={Linear Algebra Appl.},
   volume={433},
   date={2010},
   number={11-12},
   pages={1726--1759},
   issn={0024-3795},
   review={\MR{2736095}},
   doi={10.1016/j.laa.2010.06.039},
}



\end{biblist}
\end{bibdiv}
\end{document}